\documentclass{article}
\usepackage{amsfonts, amssymb, latexsym, epsfig}

\newtheorem{lemma}{Lemma}
\newtheorem{proposition}{Proposition}

\begin{document}
\title{$M$-curves of degree  $9$ or $11$ with one unique non-empty oval}
\author{S\'everine Fiedler-Le Touz\'e}
\maketitle
\begin{abstract}
In this note, we consider $M$-curves of odd degree $m$ with real scheme of the
form $\langle \mathcal{J} \amalg  \alpha \amalg 1
\langle \beta \rangle \rangle$. With help of complex orientations, we prove that
for $m = 9$, $\alpha \geq 2$, and for $m = 11$, $\alpha \geq 3$.
\end{abstract}
\section{Introduction}
Let $A$ be a real algebraic non-singular plane curve of degree $m$, its complex part $\mathbb{C}A \subset \mathbb{C}P^2$ is a Riemannian surface with genus $g = (m-1)(m-2)/2$ and its real part $\mathbb{R}A$ is a set of $L \leq g+1$ circles embedded in $\mathbb{R}P^2$. A circle embedded in $\mathbb{R}P^2$ is called {\em oval\/} or {\em pseudo-line\/} depending on whether it realizes the class 0 or 1 of $H_1(\mathbb{R}P^2)$. If $m$ is even, all of the circles are ovals, if
$m$ is odd, the curve has one unique pseudo-line, denoted by $\mathcal{J}$. 
An oval separates $\mathbb{R}P^2$  into a M\"obius band and a disc. The latter is called the {\em interior\/} of the oval. An oval of $\mathbb{R}A$ is {\em empty\/} if its interior contains no other oval. One calls {\em exterior oval\/} an oval that is surrounded by no other oval. 
We say that $A$ is an $M$-curve if $L = g+1$. A curve is {\em dividing\/} if its real part disconnects its complex part: $\mathbb{C}A \setminus \mathbb{R}A$ has then two homeomorphic halves that are exchanged by the complex conjugation. One can endow $\mathbb{R}A$ with a {\em complex orientation\/} induced by the orientation of one of these halves, note that the complex orientation is defined only up to complete reversion. The $M$-curves are dividing. If $A$ is dividing, then $L \equiv g + 1$ (mod 2).
Two ovals form an {\em injective pair\/} if one of them lies in the interior of the other one. One can provide all the injective pairs of $\mathbb{R}A$ with a sign as follows: such a pair is {\em positive\/} if and only if the orientations of its two ovals induce an orientation of the annulus that they bound in $\mathbb{R}P^2$. Let $\Pi_+$ and $\Pi_-$ be the numbers of positive
and negative injective pairs of $A$. If $A$ has odd degree, each
oval of $\mathbb{R}A$ can be endowed with a sign: given an oval $O$
of $\mathbb{R}A$, consider the M\"obius band $\mathcal{M}$
obtained by cutting away the interior of $O$ from $\mathbb{R}P^2$.
The classes $[O]$ and $[2\mathcal{J}]$ of $H_1(\mathcal{M})$ either
coincide or are opposite. In the first case, we say that $O$ is
{\em negative\/}; otherwise $O$ is {\em positive\/}. Let $\Lambda_+$ and $\Lambda_-$ be respectively the numbers of positive and negative ovals of
$\mathbb{R}A$.
Let us call the isotopy type of  $\mathbb{R}A \subset \mathbb{R}P^2$ the
{\em real scheme\/} of $A$. We consider here curves of odd degree with one unique non-empty oval, their real schemes are denoted by $\langle \mathcal{J} \amalg \alpha \amalg 1 \langle \beta \rangle \rangle$.
If $A$ is dividing, its {\em complex scheme\/} is obtained by
enriching the real scheme with the complex orientation. The 
complex schemes of our curves will be encoded by $\langle \mathcal{J} \amalg  \alpha_+ \amalg \alpha_- \amalg 1_{\epsilon} \langle \beta_+ \amalg \beta_- \rangle \rangle$ where $\epsilon \in \{ +,- \}$ is the sign of the non-empty oval,
$\alpha_+, \alpha_-$ are the numbers of positive and negative ovals
among the $\alpha$ exterior empty ovals, $\beta_+, \beta_-$ are the numbers of positive and negative ovals among the $\beta$ interior ovals. 
Let $A$ be a dividing curve of degree $m = 2k+1$, endowed with a complex orientation.
\begin{description}
\item[Rokhlin-Mishachev formula:]
{\em If\/} $m = 2k + 1$, {\em then\/}
\begin{displaymath}
2(\Pi_+ - \Pi_-) + (\Lambda_+ - \Lambda_-) = L - 1 - k(k+1).
\end{displaymath}

For the $M$-curves, one has thus $2(\Pi_+ - \Pi_-) + (\Lambda_+ - \Lambda_-)= k^2 - 2k$.

\item[Fiedler theorem:]
{\em Let $\mathcal{L}_t = \{L_t, t \in [0,1]\}$ be a pencil of real lines
based at a point $P$ of $\mathbb{R}P^2$. Consider two lines $L_{t_1}$
and $L_{t_2}$ of $\mathcal{L}_t$, which are tangent to $\mathbb{R}A$ at
two points $P_1$ and $P_2$, such that $P_1$ and $P_2$ are related by a pair
of conjugated imaginary arcs in $\mathbb{C}A \cap (\bigcup \mathbb{C}L_t)$.

Orient $L_{t_1}$ coherently to $\mathbb{R}A$ in $P_1$, and transport
this orientation through $\mathcal{L}_t$ to $L_{t_2}$.
Then this orientation of $L_{t_2}$ is compatible to that of $\mathbb{R}A$
in $P_2$.\/}
\end{description}

\begin{figure}[htbp]
\centering
\includegraphics{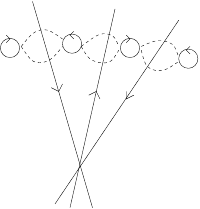}
\caption{\label{chaine} Fiedler chain of ovals}
\end{figure}

\begin{figure}[htbp]
\centering
\includegraphics{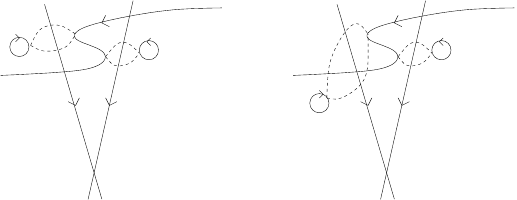}
\caption{\label{zigzag} Folds}
\end{figure}

Later on, we shall consider a complete pencil $\mathcal{L}_t$, the parameter $t$ describes a circle $\mathbb{R}P^1$ instead of an interval. 
A pencil turning in some direction gives rise to two types of tangency points, depending on whether the number of real intersection points decreases or increases by two. A sequence of ovals that are connected one to the next by pairs of conjugated imaginary arcs will be called a {\em Fiedler chain\/}, see Figure~\ref{chaine}.  
To form a chain, one allows actually two consecutive ovals to be
connected indirectly via a {\em fold\/} due to some component of the curve as shown in Figure~\ref{zigzag}. A fold is a pair of tangency points of each type on the component. To apply the theorem as a restriction tool, we try to distribute the ovals of the curve under consideration in chains, in which the ovals have alternating orientations with respect to some pencil of lines. We can safely ignore the folds.


Let $C_m$ be an $M$-curve of degree $m$. 
Given an empty oval $X$ of $C_m$, we often have to consider one point chosen in the interior of $X$. For simplicity, we call this point also $X$. 
Let $m$ be odd, for any two empty ovals $X, Y$, we shall denote by $[XY]$ and $[XY]'$ the segments of the line $XY$ that cut the pseudo-line $\mathcal{J}$ an even and an odd number of times respectively. We say that $[XY]$ is the {\em principal segment\/} determined by $X, Y$.
Let $X, Y$ and $Z$ be three interior ovals of $C_m$. The corresponding three points give rise to four triangles in $\mathbb{R}P^2$. We denote by $XYZ$ the trianlge bounded by $[XY]$, $[XZ]$, $[YZ]$ and call it the {\em principal triangle\/}.
The complete pencil of lines based at $X$ is divided in two portions by the lines $XY$ and $XZ$.
We say that the portion sweeping out $[XY]'$ has a {\em $\mathcal{J}$-jump\/}.  
For example, assume that the folds in Figure~\ref{zigzag} are on $\mathcal{J}$. The pencil of lines between the two ovals has a $\mathcal{J}$-jump in the left picture, and none in the right picture.

\section{Results}

Let $C_m$ be a curve of odd degree $m$ with one unique non-empty oval, consider the union of all principal triangles whose vertices are interior ovals.
If this union is a disc bounded by a polygon, we say that it is the {\em convex hull\/} of the interior ovals, and call it $\Delta$. It was first observed by S.~Orevkov that for $m \geq 13$, the interior ovals may have no convex hull.  

\begin{lemma}
If $m \leq 11$, then the interior ovals have a convex hull $\Delta$.
\end{lemma}

\begin{figure}[htbp]
\centering
\includegraphics{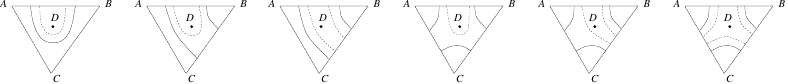}
\caption{\label{treize} The dotted arcs in the triangles are pieces of $O$, the plain arcs are pieces of $\mathcal{J}$}
\end{figure}

{\em Proof:\/}
Assume there exist four interior ovals $A, B, C, D$, such that $D$ lies in the principal triangle $ABC$. One has either $ABC = ABD \cup BCD \cup ACD$, or
$\mathcal{J}$ cuts an odd number of times each of the three segments joining $D$ to $A, B$ and $C$ in $ABC$. Assume the latter case occurs. 
Up to the action of $S_3$ on $A, B, C$, there exist a piece of $\mathcal{J}$ and a piece of the non-empty oval $O$ in $ABC$ realizing one of the six positions displayed in Figure~\ref{treize}. Note that as $A, B$ are interior ovals, the segment $[AB]$ has two supplementary intersection points with $O$. The segment $[AB]'$ cuts $\mathcal{J}$ at least once, and $O$ at least twice.
The line $AB$ cuts $C_m$ at $13$ points or more, so $m \geq 13$. For $m \leq 11$, one constructs $\Delta$ inductively, adding one oval after the other. Start with the principal triangle $\Delta_3$ determined by three ovals. Assume one has obtained with $n$ ovals a disc $\Delta_n$ bounded by a polygon. Add a new oval $X$, if $X$ is in $\Delta_n$, then $\Delta_{n+1} = \Delta_n$. Otherwise, consider the pencil of lines $\mathcal{F}_X$ sweeping out $\Delta_n$, the extreme ovals $Y, Z$ met by the pencil are vertices of $\Delta_n$. One has $\Delta_{n+1} = \Delta_n \cup XYZ$. $\Box$

\begin{lemma}
Let $C_m$ be a dividing curve of odd degree $m = 2k+1$ with real scheme $\langle \mathcal{J} \amalg  \alpha \amalg 1 \langle \beta \rangle \rangle$. If the interior ovals have a convex hull, then:
\[
1 - k - \alpha \leq \Pi_+ - \Pi_- \leq k - 1 + \alpha - \epsilon
\]
where $\epsilon = 0$ if $\alpha = 0$, $\epsilon = 1$ otherwise.
\end{lemma}

\begin{lemma}
Let $C_m$ be an $M$-curve of odd degree $m \geq 7$ with one unique non-empty oval and at least one exterior empty oval. 
\begin{enumerate}
\item
If $O$ is negative, then: $\alpha_+ - \alpha_- \equiv (k - 1)^2$ (mod $3$).
\item
Assume that the interior ovals have a convex hull.
If $O$ is positive, then:
\[
k^2 - 3k + 1 \leq 2\alpha_+ \leq 2\alpha
\]
If $O$ is negative, then:
\[
k^2 - 5k + 7 \leq 2\alpha_+ + 2\alpha \leq 4\alpha
\]
In both cases, one has: 
\[
\alpha \geq (k^2 - 5k + 7)/4
\]
\end{enumerate}
\end{lemma}

{\em Remark:\/} In degree $7$, there exist $M$-curves realizing $14$ different real schemes: $\langle \mathcal{J} \amalg \alpha \amalg 1 \langle \beta \rangle \rangle$ 
($1 \leq \alpha \leq 13$) and $\langle \mathcal{J} \amalg 15 \rangle$. The classification of the complex schemes is presented in \cite{flt1}, it was established using the restrictions and constructions from \cite{fi}, \cite{flt1}, \cite{ko}, \cite{or1} and \cite{vi1}. 
One has actually for the $M$-curves with a non-empty oval: $\Pi_+ - \Pi_- \in \{ 0, +1, -1, 2 \}$. 

\begin{figure}[htbp]
\centering
\includegraphics{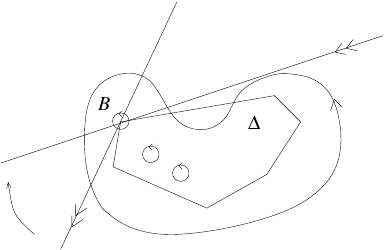}
\caption{\label{equal} The last oval of one chain is connected with $O$}
\end{figure}

{\em Proof of Lemmas~2 and 3\/:} Consider the convex hull $\Delta$ of the interior ovals. Let $B$ be an oval placed at a vertex. The pencil $\mathcal{F}_B$ sweeping out $\Delta$ sweeps out all of the interior ovals and possibly some of the $\alpha$ exterior ovals, giving rise to $k-2$ chains. We consider here not the complete chains, but the longest subchains whose extremities are interior ovals. 
There are no $\mathcal{J}$-jumps between the interior ovals, two consecutive interior ovals in a chain have thus alternating orientations with respect to $\mathcal{J}$. 
A chain has $j$ {\em jumps\/} if it is formed of $2j+1$ subchains with alternatively interior and exterior ovals. A chain with $j$ jumps brings a contribution to $\Pi_+-\Pi_-$ whose absolute value is less or equal to $j + 1$. The total number of jumps is at most $\alpha$. Taking into account the contribution of the last oval $B$ one gets the inequality $\vert \Pi_+ - \Pi_- \vert \leq k - 1 + \alpha$; if the equality holds, then the ovals vertices of $\Delta$ have all the same orientation.
Assume there exists some exterior oval $A$ lying outside of $\Delta$, and consider the pencil of lines $\mathcal{F}_A$ sweeping out $\Delta$. There are no $\mathcal{J}$-jumps between the interior ovals, and these ovals are distributed in $k-2$ chains with at most $\alpha - 1$ jumps in total. Thus $\vert \Pi_+ - \Pi_- \vert \leq k + \alpha - 3$. If the equality $\vert \Pi_+ - \Pi_- \vert = k - 1 + \alpha$ holds, then the exterior ovals lie all inside of $\Delta$ and 
some edge of $\Delta$ intersects the non-empty oval $O$. Let $B$ be a vertex of $\Delta$, extremity of such an edge. The pencil of lines based at $B$ sweeping out $\Delta$ meets all of the other empty ovals, these ovals are distributed in $k-2$ chains whose extremal ovals are all interior ovals having the same orientation as $B$. The last oval of one chain must be connected with a point on $O$, this oval forms therefore a negative pair with $O$, see Figure~\ref{equal}. Thus, 
$\Pi_+ - \Pi_-  = 1 - k - \alpha$. This finishes the proof of Lemma~2. 

The Rokhlin-Mishachev formula may be rewritten for the $M$-curves with one unique non-empty oval:
\[
\mbox{$O$ positive:} \qquad \Pi_+ - \Pi_- + 1 + \alpha_+ - \alpha_- = k^2 - 2k
\]
\[
\mbox{$O$ negative:} \qquad 3(\Pi_+ - \Pi_-) - 1 + \alpha_+ - \alpha_- = k^2 - 2k
\]

Lemma~3 follows immediately from these formulas, combined with the inequality $\Pi_+ - \Pi_- \leq k - 2 + \alpha$. 
$\Box$


{\em Remark:\/} Consider a curve $C_m$ of arbitrary degree $m$, a real line $L$, and the line pencils with base points on $L$. The sets of complex conjugated arcs obtained fiber the curve $C_m$, one may recalculate its Euler characteristic. The non-real intersection points of $C_m$ with $L$ are meridian-like singularities (like the poles of the Earth fibered by the meridians), the real intersection points are parallel-like singularities (like the poles of the Earth fibered by the parallels). Each of these singularities brings a contribution of $+1$. The real inflection points of $C_m$ and the tangencies at non real-points of $C_m$ correspond to the bifurcations, they bring each a contribution of $-1$. The tangent lines are actually bitangents. Denote by $i_r$ the number of real inflection points of $C_m$, and by $t"$ the number of bitangents to pairs of complex conjugated non-real points. One has $\chi = 2 - 2g = m -  i_r - 2t"$, hence $i_r + 2t" = m(m-2)$, this is the famous Klein formula. The fact that one can recover this formula this way was first observed by T.~ Fiedler in his phD-thesis. 

\begin{proposition}
There doesn't exist $M$-curves of degree $m = 2k+1$ ($k = 3, 4$ or $5$), having one unique non-empty oval that contains all others in its interior.
Otherwise stated, the real schemes: 
$\langle \mathcal{J} \amalg 1 \langle 14 \rangle \rangle$ ($k = 3$), 
$\langle \mathcal{J} \amalg 1 \langle 27 \rangle \rangle$ ($k = 4$), 
$\langle \mathcal{J} \amalg 1 \langle 44 \rangle \rangle$ ($k = 5$) are not realizable.  
\end{proposition}

\begin{figure}[htbp]
\centering
\includegraphics{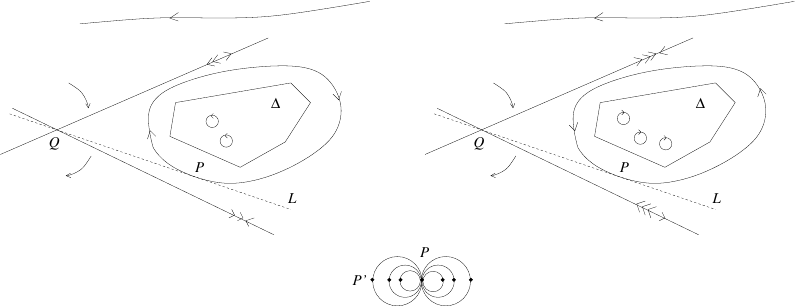}
\caption{\label{bifbis} $M$-curves $\langle \mathcal{J} \amalg 1_+ \langle 6_+ \amalg 8_- \rangle \rangle$ and $\langle \mathcal{J} \amalg 1_- \langle 15_+ \amalg 12_- \rangle \rangle$}
\end{figure}

{\em Proof\/:}
Assume there exists some curve $C_m$ contradicting the Proposition, let $O$ be its non-empty oval. Denote by $\beta_\pm$ the numbers of positive and negative interior ovals, and let $\lambda = \beta_+ - \beta_-$. For $O$ positive, $\lambda = 1 - k^2 + 2k$; for $O$ negative $\lambda = (k-1)^2/3$.
By Lemma~2, $\vert \lambda \vert = \vert \Pi_+-\Pi_- \vert \leq k-1$.
By Lemma~3, the non-empty oval can be negative only if $k \equiv 1$ (mod 3), so $k = 4$
The only possibility is a curve of degree $9$ with complex scheme $\langle \mathcal{J} \amalg 1_- \langle 15_+ \amalg 12_- \rangle \rangle$.
If $O$ is positive and $k > 3$, then $\lambda = 1 - k^2 + 2k < 1 - k < 0$, contradiction. The only possibility is a curve of degree $7$ with complex scheme
$\langle \mathcal{J} \amalg 1_+ \langle 6_+ \amalg 8_- \rangle \rangle$.
The real scheme $\langle \mathcal{J} \amalg 1 \langle 14 \rangle \rangle$ is not realizable \cite{vi1}, \cite{fi} and the complex scheme $\langle \mathcal{J} \amalg 1_- \langle 15_+ \amalg 12_- \rangle \rangle$ has been excluded more recently in \cite{go}. We give here an alternative method to exclude both complex schemes. 
Let $C_m$, $m =7$ or $9$ realize one of these complex schemes.
Let $P$ be a point of the non-empty oval $O$ such that $O$ is locally convex at $P$: if one chooses a point $P'$ close enough to $P$ on $O$, the principal segment $[PP']$ lies inside of $O$. Let $L$ be the line tangent to $O$ at $P$. We consider the 
complete pencils of lines $\mathcal{L}_{t, Q}$ based at a mobile point $Q$ that percourses the line $L$, with starting end ending position at $P$. Each of the pencils contains the line $L$, tangent to $O$ at $P$. Near the start, $Q$ is close to $P$ on $L$, another line of the pencil is tangent to $O$, at a point $P'$ close to $P$, and the two points $P, P'$ are connected by a pair of conjugated arcs forming a small circle in $\mathbb{C}A \cap (\bigcup \mathbb{C}L_{t, Q})$. As $Q$ moves away from $P$, one gets a set of growing circles, all tangent to one another at $P$. Near the end of the motion, when $Q$ moves towards $P$ from the other side of $L$, we have symmetrically a set of shrinking circles, that lie on the other side of $P$ on $\mathbb{C}A$ as shown in Figure~\ref{bifbis}.    
As $\mathbb{C}A$ is not a sphere, there must have been at some moment a bifurcation in the set of conjugated imaginary arcs, after which the two tangency points $P, P'$ with $O$ are no longer connected one to the other.
Each tangency point is then connected to an interior oval whose orientation (with respect to $\mathcal{J}$) coincides with that of $O$, this oval is the starting oval of one of the $k-1$ interior chains determined by the pencil. But for any choice of base point outside of $O$, the interior ovals are distributed in $k-1$ chains whose ovals have alternating orientations with respect to $\mathcal{J}$, each chain bringing a contribution $+1$ to $\Pi_+-\Pi_-$. The bifurcation is impossible, see 
Figure~\ref{bifbis}. (See also \cite{fi} where a similar argument is used in a different proof).
$\Box$


\begin{proposition}
There doesn't exist $M$-curves of degree $2k + 1$ with the following real schemes:
$\langle \mathcal{J} \amalg 1 \amalg 1 \langle 26 \rangle \rangle$ ($k = 4$);
$\langle \mathcal{J} \amalg 1 \amalg 1 \langle 43 \rangle \rangle$,
$\langle \mathcal{J} \amalg 2 \amalg 1 \langle 42 \rangle \rangle$ ($k = 5$).  
\end{proposition}

{\em Proof:\/}
For $k = 5$, one has by Lemma~3: $\alpha \geq 2$ if $O$ is negative, and $\alpha \geq 6$ if $O$ is positive.
Assume $O$ is negative and $\alpha = 2$. As $\alpha_+ - \alpha_- \equiv 1$ (mod $3$), one must have $\alpha_+ = 0$, $\alpha_- = -2$. 
Hence, $2\alpha_+ + 2 \alpha = 4 < k^2 - 5k + 7$, contradiction. 
For $k = 4$, one has by Lemma~2: $\alpha_+ \geq 3$ or $\alpha_+ - \alpha_- \equiv 0$ (mod $3$). Neither condition is realized if $\alpha = 1$. $\Box$

Note that for $m = 7$, the real scheme $\langle \mathcal{J} \amalg 1 \amalg 1 \langle 13 \rangle \rangle$
is realizable, with the two admissible complex schemes:
$\langle \mathcal{J} \amalg 1_+ \amalg 1_+ \langle 6_+ \amalg 7_ -\rangle \rangle$ and
$\langle \mathcal{J} \amalg 1_+ \amalg 1_- \langle 7_+ \amalg6 _- \rangle \rangle$, see \cite{vi2}, \cite{ko}, \cite{flt1}. 
For $m = 9$, the non-realizability of the schemes 
$\langle \mathcal{J} \amalg 1 \langle 27 \rangle \rangle$,
$\langle \mathcal{J} \amalg 1 \amalg 1 \langle 26 \rangle \rangle$
proved above had previously been announced in \cite{fi}, but unfortunately the proofs were never published and went lost. 
Ninth degree $M$-curves $\langle \mathcal{J} \amalg \alpha \amalg 1 \langle \beta \rangle \rangle$ with $1 \leq \beta \leq 19$, $\beta = 22, 23$ have been constructed by A.~Korchagin \cite{ko6}. 
By Lemma~3, the only admissible schemes for $\alpha = 2$ and $3$ are:
$\langle \mathcal{J} \amalg 1_+ \amalg 1_- \amalg 1_- \langle 14_+ \amalg
11 _- \rangle \rangle$, 
$\langle \mathcal{J} \amalg 3_+ \amalg 1_+ \langle 10_+ \amalg
14 _- \rangle \rangle$, 
$\langle \mathcal{J} \amalg 3_+ \amalg 1_- \langle 13_+ \amalg
11 _- \rangle \rangle$, 
$\langle \mathcal{J} \amalg 3_- \amalg 1_- \langle 14_+ \amalg
10 _- \rangle \rangle$.
Finally, note that the proofs and results presented here are valuable also for pseudo-holomorphic curves.

\vspace{2cm}

severine.fiedler@live.fr

\end{document}